\documentclass[12pt]{article}
\usepackage{latexsym}

\def\bs{\begin{subequations}}
\def\es{\end{subequations}}

\catcode`\@=11

\newtoks\@stequation
\def\subequations{\refstepcounter{equation}
  \edef\@savedequation{\the\c@equation}%
  \@stequation=\expandafter{\theequation}
  \edef\@savedtheequation{\the\@stequation}
  \edef\oldtheequation{\theequation}%
  \setcounter{equation}{0}%
  \def\theequation{\oldtheequation\alph{equation}}}

\def\endsubequations{\setcounter{equation}{\@savedequation}%
  \@stequation=\expandafter{\@savedtheequation}%
  \edef\theequation{\the\@stequation}\global\@ignoretrue}

\makeatletter
        \renewcommand{\theequation}{\thesection.\arabic{equation}}%
        \@addtoreset{equation}{section}%
\makeatother

\renewcommand{\thefootnote}{\fnsymbol{footnote}}

\begin{document}

\begin{titlepage}

April 2, 2008  added Appendix B: SU(2)

\begin{center}        \hfill   \\
            \hfill     \\
                                \hfill   \\

\vskip .25in

{\large \bf Analytic Functions of a Quaternionic Variable \\}

\vskip 0.3in

Charles Schwartz\footnote{E-mail: schwartz@physics.berkeley.edu}

\vskip 0.15in

{\em Department of Physics,
     University of California\\
     Berkeley, California 94720}
        
\end{center}

\vskip .3in

\vfill

\begin{abstract}

Here we follow the basic analysis that is common for real and 
complex variables and find how it can be applied to a quaternionic variable.
Non-commutativity of the quaternion algebra poses obstacles for the 
usual manipulations; but we show how  
many of those obstacles can be overcome. After 
a tiny bit of linear algebra we look at the beginnings of 
differential calculus. The surprising result is that the first order 
term in the expansion of $F(x+\delta)$ is a compact formula 
involving both $F^{\prime}(x)$  
and $[F(x) - F(x^{*})]/(x-x^{*})$.

\end{abstract}

\vfill

\end{titlepage}

\renewcommand{\thefootnote}{\arabic{footnote}}
\setcounter{footnote}{0}
\renewcommand{\thepage}{\arabic{page}}
\setcounter{page}{1}

\section{Introduction}
We are very familiar with functions of a real or complex variable $x$ 
which we can expand, in the mode of differential calculus, as
\begin{equation}
F(x+\delta) = F(x) + F^{\prime}(x) \delta + \frac {1}{2}F^{\prime 
\prime}(x) \delta^{2} + \ldots.\label{1}
\end{equation}

But what if we consider a quaternionic variable
\begin{equation}
x = x_{0} + i x_{1} + j x_{2} + k x_{3} ,\label{2}
\end{equation}
where those quaternions $i, j,k$ do not commute with one another. The 
small quantity $\delta$ will also involve all those quaternions. How 
then can we expect anything as neatly packaged as Equation (\ref{1})?

This is a long-standing challenge to mathematicians; and here we 
believe we have something new to offer. Previous work has attempted 
to extend the usual concepts of the derivative, $\frac{dF}{dx}$: see, 
for example, references \cite{1}, \cite{2} and \cite{3}. We shall, instead,  
focus only on the 
differential $dF$, as suggested by Eq. (\ref{1}). The method will be to start 
by examining the exponential function $e^{x}$, in Section 3, 
and then use this result to build more general functions $F(x)$, 
 in Section 4. We succeed in finding a general formula for the first 
 order term that 
is surprisingly compact, as shown in Eq. (\ref{c3a});  but the price we pay 
is that it is no longer local, involving $[F(x) - F(x^{*})]/(x-x^{*})$ along with 
$F^{\prime}(x)$.

In Section 5 we describe a geometric view of this situation;
in Section 6 we find the Leibnitz rule to be valid. Section 7 
presents the second order term in the generalized expansion (\ref{1}) 
and Section 8 contains a general discussion of these results.
In Appendix B we show how this line of analysis can be extended to 
the non-commutative variables based on the Lie algebra SU(2,C).

First, however, we make a brief visit to quaternionic linear algebra.

\section{Linear Algebra}

Suppose we have the  linear equation
\begin{equation}
ax + xb = c \label{a1}
\end{equation}
where $a,b,c$ are given quaternions and we want to solve for the 
unknown quaternion $x$. If we were dealing with real or complex 
quantities, then all the factors in Equation (\ref{a1}) would commute 
and we would write $x = (a+b)^{-1}\;c$. But that does not work for 
quaternionic numbers and variables. 

In fact, the solution to Equation(\ref{a1}) may be written as
\begin{equation}
x = d^{-1}(ac + cb^{*}), \;\;\; d=a^{2} + a(b+b^{*}) + b^{*}b 
\label{a2}
\end{equation}
where one recognizes that $[a,d]=ad-da=0$. The symbol $^{*}$ means the 
usual complex conjugation, which changes the signs of all the 
imaginaries. I would imagine that 
this solution is not new but I do not know where to find that out.

There is an alternative way of writing the solution:
\begin{equation}
x = (a^{*}c + cb)h^{-1},\label{a3}
\end{equation}
which I will leave for the reader to explore.

\section{Exponential Function}

For a general quaternionic variable x, we can define the exponential 
function in the usual way
\begin{equation}
e^{x} = \lim_{N\rightarrow \infty} (1+\frac{x}{N})^{N} \label{b1}
\end{equation}
and this leads us to the expansion
\begin{equation}
e^{(x+\delta)} = e^{x}[1 + \int_{0}^{1}ds\; e^{-sx}\; \delta\;e^{sx} 
+ O(\delta^{2})].\label{b2}
\end{equation}
This formula is well known to some people (I believe it is often 
credited to Richard Feynman but I do not know a reference);  
 a derivation of it is given in Appendix A.

If we now separate x into real and imaginary parts as
$x = x_{0} + r u_{x}$, where $(u_{x})^{2} = -1$, and, furthermore, 
use the expansion 
\begin{equation}
e^{x} = e^{x_{0}}( cos\;r + u_{x}\;sin\;r) ,\label{b3}
\end{equation}
then we can write the expansion as
\begin{equation}
e^{(x+\delta)} = e^{x} + e^{x}[a \;\delta + b\; [u_{x},\delta] +
c\; u_{x}\; \delta\; 
u_{x} + O(\delta^{2})]\label{b4}
\end{equation}
where the real quantities $a,b,c$ are given by
\begin{equation}
a = \frac{1}{2}\;(1+ \frac{sin2r}{2r}), \;\;\;
b = \frac{1}{2}\;\frac{cos2r -1}{2r}, \;\;\; 
c = \frac{1}{2}\;(-1 + \frac{sin2r}{2r}).\label{b5}
\end{equation}
In the special case where $u_{x}$ commutes with 
$\delta$, then this reduces to the familiar formula for the 
exponential of a complex variable.

We can now go on to the logarithm function, defined as the inverse of the 
exponential:
\begin{equation}
e^{x} = y, \;\;\; x = ln \;y,\;\;\; e^{(x+\delta)} = y + \Delta, \;\;\; 
x+\delta = ln (y+\Delta).\label{b6}
\end{equation}

Using the previous results, we eventually arrive at the expansion 
\begin{eqnarray}
ln(y+\Delta) = ln \;y + A y^{-1}\Delta + B [u_{x},y^{-1}\Delta] + C 
u_{x}\;y^{-1}\Delta\;u_{x} + 
O(\Delta^{2}),\nonumber \\
A =\frac{1}{2}(r\; cot r +1), \;\;\; B = \frac{r}{2}, \;\;\; 
C=\frac{1}{2}(r\; cot r -1).\label{b7}
\end{eqnarray}

\section{General Function F(x)}

For a general function $F(x)$ of a quaternionic variable $x = 
x_{0} + r u_{x}$, we 
start by assuming a representation as a Laplace transform:
\begin{equation}
F(x) = \int dp\; f(p)\; e^{px}\label{c1}
\end{equation}
where $p$ is a real variable.  We then calculate an expansion for 
$F(x+\delta)$ using the formula derived above for the exponential 
function; and this leads to the following general formula
\begin{eqnarray}
F(x+\delta) = F(x) + F^{\prime}(x)\;\delta + \frac{1}{4r}\;(-F(x) + 
F(x^{*}))\;[u_{x},\delta] +\nonumber \\  \frac{1}{4} 
F^{\prime}(x) \;[u_{x},[u_{x},\delta]] + O(\delta^{2}),\label{c2}
\end{eqnarray}
where $F^{\prime}(x)$ is the usual derivative of the function $F(x)$ 
calculated as if $x$ were a real variable. It is surprising how 
simple this formula appears.  

 Another form of this formula (\ref{c2}) is
 \begin{eqnarray}
	F(x+\delta) - F(x) = F^{\prime}(x)\;\delta_{1} +
  (F(x) - F(x^{*}))\;(x - x^{*})^{-1} \; \delta_{2} + O(\delta^{2}), 
  \label{c3a} \\
\delta_{1} =   \frac{1}{2}\;(\delta - u_{x} \;\delta\;u_{x}) , 
\;\;\;\;\;\;\;\;  
\delta_{2} =  \frac{1}{2}\;(\delta + 
	u_{x}\delta u_{x}) ,\label{c3b}
\end {eqnarray}
which will be discussed further in the next Section. 

If one considers the function $F(x) = x^{n}$, then the above formula,  
(\ref{c2}),  is correct, as may be shown by induction. Thus it is also true 
for any power series, $F(x) = \sum_{n}\;c_{n}\;x^{n}$.

It is an interesting exercize to show that the special formula 
(\ref{b7}) for the logarithm 
function is in agreement with the general formula (\ref{c2}).

\section{Geometric View}
The variable x, as written in Eq. (\ref{2}), may be viewed as a point 
in a four-dimensional space composed of the real line, for $x_{0}$, 
and a three-dimensional Euclidean space, for the vector $\textbf{x} = 
(x_{1},x_{2},x_{3})$. In the decomposition $x = x_{0} + r u_{x}$, we 
see $r$ as the length of the vector $\textbf{x} = r \hat{x}$; and the 
unit imaginary quaternion $u_{x}$ is the dot product of the unit 
vector $\hat{x}$ with the vector $(i,j,k)$.

Now we want to recognize that the quaternion $\delta$, which is used 
to displace the variable $x$ when we write $F(x+\delta)$, can be 
given a different decomposition in that four-dimensional space: as 
defined in the previous section, Eq. (\ref{c3b}), 
$\delta = \delta_{1} + \delta_{2}$. We can recognize that 
$\delta_{1}$ is in the two-dimensional space composed of the real 
line and the radial direction along the vector $\textbf{x}$; and 
$\delta_{2}$ is in the two-dimensional space that is orthogonal to 
the direction of $\textbf{x}$.  One can say, geometrically, that 
$\delta_{1}$ is parallel to $x$ (perhaps writing it as 
$\delta_{\parallel}$) and $\delta_{2}$ is perpendicular 
($\delta_{\perp}$). 
Algebraically, this is represented by 
the  relations,
\begin{equation}
\delta_{1}\;x = x\;\delta_{1}, \;\;\;\;\; \delta_{2}\;x= 
x^{*}\;\delta_{2}.\label{c4}
\end{equation}

This decomposition of $\delta$ is a local process. It  
varies from one place $x$ to another, rather like the familiar unit 
vectors in polar coordinates of two-dimensional Euclidean space.

The authors of reference \cite{1} have taken a similar approach, 
introducing a local unit imaginary (which they call \emph{iota}) that is 
the same as what we have defined as $u_{x}$. However, they limit 
their differentiations to displacements that are restricted to this 
two-dimensional space: they allow only what we call 
$\delta_{1}$ without any of $\delta_{2}$. In that way they merely 
reproduce what is known about ordinary complex variables.

\section{Leibnitz' Rule}
Let us now define the first-order differential operator ${\cal{D}}$, from 
Eq.(\ref{c3a}), as
\begin{equation}
F(x+\delta) = F(x) + {\cal{D}}\;F(x) + O(\delta^{2})\label{c5}
\end{equation}
with
\begin{equation}
{\cal{D}}\;F(x) = F^{\prime}(x)\;\delta_{1} + 
(F(x) - F(x^{*}))\;(x - x^{*})^{-1} \; \delta_{2}.\label{c6}
\end{equation}
Now we consider ${\cal{D}} (F(x)G(x))$. It starts off easy:
\begin{equation}
{\cal{D}} (FG) = (F(x) + {\cal{D}} F)(G(x) + {\cal{D}} G) - FG = F {\cal{D}} G + 
({\cal{D}} F)G;\label{c7}
\end{equation}
but the next step is more interesting, substitution (\ref{c6}) into 
(\ref{c7}):
\begin{eqnarray}
{\cal{D}} (FG) =&& F [G^{\prime}(x)\;\delta_{1} + 
(G(x) - G(x^{*}))\;(x - x^{*})^{-1}\delta_{2}] +\nonumber \\  
 &&[F^{\prime}(x)\;\delta_{1} + 
(F(x) - F(x^{*}))\;(x - x^{*})^{-1}\delta_{2}]G.\label{c8}
\end{eqnarray}
The question is, How do we move the quaternions $\delta$ to the 
right, past the function $G$ in the second term of Eq. 
(\ref{c8})? 
The answer is given by (\ref{c4}); and we finally get,
\begin{equation}
{\cal{D}}(FG) = (F\;G^{\prime} + F^{\prime}\;G)\delta_{1} + 
[F(x)G(x) - F(x^{*})G(x^{*})](x-x^{*})^{-1}\;\delta_{2},
\end{equation}
which corroborates Leibnitz' rule for this differential operator.

\section{Second Order Terms}

Let's return to the exponential function (\ref{b1}) and proceed with 
 the expansion,
\begin{equation}
e^{(x+\delta)} = e^{x}[1+\int_{0}^{1}ds\; e^{-sx}\delta e^{sx} +
\int_{0}^{1}dt \int_{0} ^{1-t}ds \; e^{-(s+t)x}\delta e^{tx}\delta 
e^{sx} + O(\delta^{3})].
\end{equation}

The better approach is to combine the exponential function and the 
Laplace transform from the beginning. Writing $F(x+\delta) = F(x) + 
F^{(1)} + F^{(2)} + \ldots$, we will first re-do the calculation of 
$F^{(1)}$ to see how easily it goes with the decomposition $\delta = 
\delta_{1}+\delta_{2}$ and the algebra (\ref{c4}).
\begin{eqnarray}
F^{(1)} = \int dp f(p) \;p\int _{0}^{1}ds\; e^{(1-s)px}(\delta_{1}+ 
\delta_{2})\;e^{spx} = \\ 
\int dp f(p)\;p \;e^{px}\;\int _{0}^{1}ds[\delta_{1} 
+e^{sp(x^{*}-x)}\;\delta_{2}]
\end{eqnarray}
and this leads immediately to the result Eq. (\ref{c3a}).

Now we look at
\begin{equation}
F^{(2)} =\int dp f(p)\; e^{px} p^{2}\;
\int_{0}^{1}dt \int_{0} ^{1-t}ds \; e^{-(s+t)px}\delta e^{tpx}\delta 
e^{spx}.
\end{equation}
Again, we decompose $\delta$ and after a bit more work arrive at the 
result for the second order term,
\begin{eqnarray}
F^{(2)} = \frac{1}{2} F^{\prime \prime}(x)\;\delta_{1}^{2} + 
(F(x)-F(x^{*}))\;(x-x^{*})^{-2}\;(\delta_{2}\delta_{1} - \delta 
\delta_{2}) + \nonumber \\ 
F^{\prime}(x)\;(x-x^{*})^{-1}\;\delta\delta_{2} +
F^{\prime}(x^{*})\;(x^{*}-x)^{-1}\;\delta_{2}\delta_{1}.
\end{eqnarray}
\vskip 1cm

\section{Discussion}

One may ask what restrictions there are on the functions $F(x)$ 
considered above. At first, one would say that they should be real 
analytic 
functions; having terms like $xax$ where a is a general quaternion 
would certainly cause trouble. One can extend this condition slightly 
by allowing $F(x)$ (but not the function $G(x)$ in Section 6) 
to be a real function with arbitrary quaternions 
multiplying from the left. That is, the power series form $F = 
\sum_{n}\;c_{n}\;x^{n}$ could have arbitrary numbers $c_{n}$.

This bias to the left-hand side can be reversed if we change the 
original steps (\ref{b2}), setting $s \rightarrow 1-s$,  and 
(\ref{c1}), putting $f(p)$ on the right-hand side.

It is noteworthy that our differential operators are no longer local: 
they involve $F(x^{*})$ along with $F(x)$. Nevertheless, it is 
surprising how simple and how general the results obtained here are.

The Taylor series we have discussed above are expansions about the 
origin $x=0$. In the usual complex analysis such power series may be 
about any fixed point $x=x_{f}$; but such a quaternion constant put in 
the middle of our expressions would appear to cause trouble. That can 
be rectified by defining the new variable $y = x - x_{f} = y_{0} + 
r_{y}\;u_{y}$ and then using this new imaginary $u_{y}$ to separate 
the displacement $\delta = \delta_{1}+\delta_{2}$.

One may also ask if this general method may be applied to some other 
kind of non-commuting algebra beyond the quaternions.  I believe that 
something very similar can be done starting with a Clifford algebra.
Another example is given in Appendix B.

\vskip 0.5cm
\begin{center} \textbf{ACKNOWLEDGMENT} \end{center}

 I am grateful to J. Wolf for some helpful conversation.

\vskip 0.5cm
\setcounter{equation}{0}
\def\theequation{A.\arabic{equation}}
\boldmath
\noindent{\bf Appendix A}
\unboldmath
\vskip 0.5cm

Here we give a derivation of the formula (\ref{b2}) for any 
non-commuting quantities $x$ and $\delta$.
\begin{eqnarray}
e^{(x+\delta)} = \lim_{N\rightarrow \infty} [1+\frac{x}{N} + 
\frac{\delta}{N}]^{N} = \\ 
\lim_{N\rightarrow \infty}\{[1+\frac{x}{N}]^{N} +
\sum_{m=0}^{N-m-1}\;[1+\frac{x}{N}]^{N-m-1} \;\frac{\delta}{N}\; [1+\frac{x}{N}]^{m} + 
O(\delta^{2})\}.
\end{eqnarray}
In taking the limit $N \rightarrow \infty$, we convert the sum over 
$m$ to an integral over $s = \frac{m}{N}$ and this yields
\begin{equation}
e^{(x+\delta)} = e^{x} + \int_{0}^{1}\; ds\; e^{(1-s)x}\;\delta \; 
e^{sx} + O(\delta^{2}).
\end{equation}

\vskip 0.5cm
\setcounter{equation}{0}
\def\theequation{B.\arabic{equation}}
\boldmath
\noindent{\bf Appendix B}
\unboldmath
\vskip 0.5cm

Here we shall extend the general method used above for a quaternionic 
variable to something built on a more general Lie Algebra - 
specifically SU(2).

Here is the Lie algebra:
\begin{equation}
[J_{1},J_{2}]=J_{3}, \;\;\; [J_{2},J_{3}]=J_{1},\;\;\; 
[J_{3},J_{1}]=J_{2},\label{B1}
\end{equation}
where the three  $J$'s are understood to be matrices over 
the complex numbers.  In particular we shall use the relations
\begin{equation}
e^{\theta J_{3}}\;J_{1}\;e^{-\theta J_{3}} = J_{1}\;cos\theta + 
J_{2}\;sin\theta, \;\;\;\;\; 
e^{\theta J_{3}}\;J_{2}\;e^{-\theta J_{3}} = J_{2}\;cos\theta - 
J_{1}\;sin\theta,\label{B2}
\end{equation}
which follow from (\ref{B1}).

The new variable $x$ is to be constructed with four real parameters as
\begin{equation}
x = x_{0}\;I + x_{1}\;J_{1} + x_{2}\;J_{2}+ x_{3}\;J_{3} \label{B3}
\end{equation}
and we want to expand $F(x+\delta) = F(x) + F^{(1)} + O(\delta^{2})$,
where $\delta$ is a small quantity in that same space of matrices as $x$.
Our first step is to define a local coordinate system at the given 
point $x$. By a suitable linear transformation (rotation) of the Lie 
algebra we make the coordinate $x$ appear as
\begin{equation}
x = x_{0}\;I + r J_{3}\label{B4}
\end{equation}
where we recognize that $r^{2} = x_{1}^{2}+x_{2}^{2}+x_{3}^{2}$.

We can now separate the displacement $\delta = \delta_{\parallel} + 
\delta_{\perp}$ as follows.
\begin{equation}
\delta_{\parallel} = \delta_{0}\;I + \delta_{3}\;J_{3}, \;\;\;\;\;
\delta_{\perp} = \delta_{1}\;J_{1} + \delta_{2}\;J_{2}.\label{B5}
\end{equation}

Now we are ready to study the first order term in the expansion, 
again using the representation of $F(x)$ in terms of the exponential 
function.
\begin{equation}
F^{(1)} = \int dp\;f(p)\;p\;e^{px}\;\int_{0}^{1}ds\; e^{-spx}\;\delta 
\;e^{spx}.\label{B7}
\end{equation}
Since $\delta_{\parallel}$ commutes with $x$, the first part of this 
is simply $F^{\prime}(x)\;\delta_{\parallel}$.  For the part with 
$\delta_{\perp}$ we use the formulas (\ref{B2}), where $\theta$ is 
replaces by $-spr$.  The integrals over $s$ are trivial and we 
merely write $sin(pr)$ and $cos(pr)$ in terms of $e^{\pm ipr}$ to get 
our final result.
\begin{eqnarray}
F(x+\delta) = F(x) + F^{\prime}(x)\delta_{\parallel} + 
\{F(x+ir)-F(x-ir)\}\;\frac{1}{2ir}\;\delta_{\perp} + \nonumber \\
\{F(x+ir) + F(x-ir) - 2F(x)\} \;\frac{1}{2r}\; [J_{3},\delta_{\perp}] +\;\;
 O(\delta^{2}).\;\;\;\;\;\label{B8}
\end{eqnarray}

It should be noted that the $\delta$-related factors in Eq. 
(\ref{B8}) can be written in the following way:
\begin{eqnarray}
[J_{3},\delta_{\perp}] = \frac{1}{r}\;[x,\delta] , \label{B9}\\ 
\delta_{\perp} = - \frac{1}{r^{2}}\;[x,[x,\delta]],\label{B10} \\
\delta_{\parallel} = \delta - \delta_{\perp}\label{B11}.
\end{eqnarray}
This means that we do not have to carry out the ''rotation'' that 
gave us Eq. (\ref{B4}) explicitely; the talk about choosing a local 
coordinate system is merely rhetorical. 

I expect that this method can be extended to other Lie algebras, with 
the quantity $\delta_{\perp}$ subdivided into distinct portions according to 
the roots of the particular algebra. The system of Eqs. 
(\ref{B10}), (\ref{B11}) would be adapted to make those separations, 
using the known values of the roots; and those root values would also 
appear in the final generalization of Eq. (\ref{B8}).

\end{document}